\documentclass{amsart}

\usepackage[all]{xy}
\usepackage{amsmath}
\usepackage{amssymb}

\newtheorem{theorem}{Theorem}[section]
\newtheorem{lemma}[theorem]{Lemma}

\theoremstyle{definition}

\newtheorem{proposition}[theorem]{Proposition}

\theoremstyle{remark}
\newtheorem{remark}[theorem]{Remark}

\numberwithin{equation}{section}

\begin{document}
\author{Yanlong Hao}
\address{School of Mathematical Sciences,
         Tianjin 300071, P.R.China}
\email{haoyanlong13@mail.nankai.edu.cn}
\author{Xiugui $\rm{Liu}^*$}
\address{School of Mathematical Sciences and LPMC,
         Tianjin 300071, P.R.China}
\email{xgliu@nankai.edu.cn}
\thanks{* The author was supported in part by the National Natural Science Foundation of China (No. 11171161), Program for New Century Excellent Talents in University (No. NCET-08-0288), and the Scientific Research Foundation for the Returned Overseas Chinese Scholars, State Education Ministry.}
\author{Qianwen Sun}
\address{School of Mathematical Sciences,
         Tianjin 300071, P.R.China}
\email{qwsun13@mail.nankai.edu.cn}
\date{}
\title{On the homotopy fixed point sets of spheres actions on rational complexes}
\maketitle
$\textbf{Abstract}$: In this paper, we describe the homotopy type of the homotopy fixed point sets of $S^3$-actions on rational spheres and complex projective spaces, and provide some properties of $S^1$-actions on a general rational complex.
\section{Introduction}
An action of a group $G$ on a space $M$ gives rise to two natural spaces, the fixed point set $M^{G}$ and the homotopy fixed point set $M^{hG}$.
It is crucially important that there is an injection:
$$k:M^G\longrightarrow M^{hG}.$$
Indeed, one version of the \emph{generalized Sullivan conjecture} asserts that, when $G$ is a finite $p$-group, and $M$ is a $G$-CW-complex, then the $p$-completion of $k$ is a homotopy equivalence. This conjecture was proved in the case when $M$ is a finite complex by Miller \cite{H}.

For a finite group $G$, the rational homotopy theory of $M^{hG}$ has been studied by Goyo \cite{J}.

In \cite{U-Y-S, U-Y-A}, the authors studied the homotopy type of $M^{hG}$ for a compact Lie group $G$ with particular emphasis when $G$ is the circle.

From now on, and unless explicitly stated otherwise, $G$ will denote a compact connected Lie group and by a topological $G$-space we mean a nilpotent $G$-space with the homotopy type of a CW-complex of finite type and $M^{G}\neq \emptyset$. Then the action of $G$ on $M$ induces an action of $G$ on $M_{\mathbb{Q}}$.

We then start by setting a sufficiently general context in which ${M_{\mathbb{Q}}}^{hG}$ has the homotopy type of a nilpotent  CW-complex. Identifying the homotopy fixed point set with the space sec($\xi$) of sections of the corresponding Borel fibration
$$\xi:M\rightarrow M_{hG}\rightarrow BG,$$
we have that if $\pi_{>n}(M)$ are torsion groups for a certain $n>1$, then ${M_{\mathbb{Q}}}^{hG}$ is a rational nilpotent complex with the homotopy type of a CW-complex \cite{U-Y-S}.

In this paper, we explicitly describe the rational homotopy type of the homotopy fixed point sets of certain $S^3$ actions.
\begin{theorem}\label{2}
$\\$
(1) When $n$ is odd, ${S^{n}_{\mathbb{Q}}}^{hS^3}$ has the rational homotopy type of products of odd dimensional spheres, precisely, we have
$${S^{n}_{\mathbb{Q}}}^{hS^3}\simeq_{\mathbb{Q}} S^a\times S^{a+4}\times\cdots \times S^n,$$
where
$$a=\left\{
\begin{array}{ll}
1,&n=4k+1,\\
3,&n=4k+3.
\end{array}
\right.  $$
(2) If $n=4k$, ${S^{n}_{\mathbb{Q}}}^{hS^3}$ is either path connected, and of the rational homotopy type of $S^3\times K_{k}$,
where $K_{k}$ has the minimal Sullivan model
$$(\Lambda ((x_{s})_{1\leq s\leq k},(y_r)_{2\leq r\leq 2k}),d)$$
with $|x_s|=4s$, $|y_{r}|=4r-1$, $dx_s=0$ $(1\leq s\leq k)$,
$\displaystyle{dy_r=\sum_{s+t=r}x_sx_t}$ $(2\leq r\leq 2k)$, or else, it has $2$ components, each of them has the rational homotopy type of
$$S^{4k+3}\times S^{4k+7}\times\cdots \times S^{8k-1},$$.\\
(3) If $n=4k+2$, ${S^{n}_{\mathbb{Q}}}^{hS^3}$ is path connected, and of the rational homotopy type of $S^3\times S^7 \times T_{k}$, where $T_{k}$ has the minimal Sullivan model
$$(\Lambda ((x_{s})_{1\leq s\leq k},(y_r)_{3\leq r\leq 2k+1}),d)$$
with $|x_s|=4s+2$, $|y_{r}|=4r-1$, $dx_s=0$ $(1\leq s\leq k)$,
$\displaystyle{dy_r=\sum_{s+t=r-1}x_sx_t}$ $(3\leq r\leq 2k+1).$
\end{theorem}
\begin{theorem}\label{4}
$\\$
(1) If $n$ is odd, ${\mathbb{C}P^{n}_{\mathbb{Q}}}^{hS^3}$ is path connected, and has the rational homotopy type of one of the following spaces:
$$\mathbb{C}P^1\times S^{7}\times S^{11}\times\cdots \times S^{2n+1},$$
$$S^3\times \mathbb{C}P^3\times S^{11} \times \cdots \times S^{2n+1},$$
$$S^3\times S^7\times \mathbb{C}P^{5} \times \cdots \times S^{2n+1},$$
$$\cdots,$$
$$S^3\times S^7\times \cdots \times S^{2n-3} \times\mathbb{C}P^{n}.$$
(2) If $n$ is even, ${\mathbb{C}P^{n}_{\mathbb{Q}}}^{hS^3}$ is path connected, and has the rational homotopy type of one of the following spaces:
$$*\times S^{5}\times S^{9}\times\cdots \times S^{2n+1},$$
$$S^1\times \mathbb{C}P^2\times S^{9} \times \cdots \times S^{2n+1},$$
$$S^1\times S^5\times \mathbb{C}P^{4} \times \cdots \times S^{2n+1},$$
$$\cdots,$$
$$S^1\times S^5\times \cdots \times S^{2n-3} \times \mathbb{C}P^{n}.$$
\end{theorem}

In \cite[Corollary 2]{U-Y-S}, they give a criterion of an elliptic $S^{1}$-space. We first show that the condition $M$ is a finite complex is necessary by the following example: there is a nilpotent $S^1$-complex $M$ which is not an elliptic space, such that each component of ${M_{\mathbb{Q}}}^{hS^1}$ is elliptic. We also observe that an $S^1$-finite nilpotent complex $M$ is elliptic if and only if one of the component of ${M_{\mathbb{Q}}}^{hS^1}$ is elliptic, complementing
the mentioned result.

Finally, we show that the injection $k$ is generally not a rational homotopy equivalence.
\begin{theorem}\label{5}
For an $S^1$-complex $M$ which is simply connected with
$${\rm{dim}}\pi_*(M)\otimes \mathbb{Q}<\infty.$$
Then
$$k:M_{\mathbb{Q}}^{S^1}\hookrightarrow M_{\mathbb{Q}}^{hS^1}.$$
is a rational homotopy equivalence if and only if $M$ is rational homotopy equivalent to a product of $CP^{\infty}$.
\end{theorem}

In the next section we prove Theorem \ref{2} and \ref{4}. In section 3 we prove Theorem \ref{5}.,
\section{$S^3$-rational spheres and complex projective spaces}
Our results heavily depend on known facts and techniques arising from rational homotopy theory. All of them can be found with all details in \cite{Y-S}. We simply remark a few facts.

We recall that when $M$ is path connected, the Sullivan model of $M$ is a quasi-isomorphism:
$$m:(\Lambda V_{M},d)\rightarrow A_{PL}(M),$$
where $(\Lambda V_{M},d)$ is a Sullivan algebra.

We also recall that a space $M$ is elliptic if both $H^*(M;\mathbb{Q})$ and $\pi_{*}(M)\otimes\mathbb{Q}$ are finite dimensional vector spaces over $\mathbb{Q}$.

For a $G$-space $M$, we have the corresponding Borel fibration:
$$\xi:M\rightarrow M_{hG}\rightarrow BG,$$
where $M_{hG}=(M\times EG)/G$. It is a classical fact that the homotopy fixed point set
$$M^{hG}={\rm{map}}_{G}(EG,M)$$
is homotopy equivalent to the section space Sec($\xi$) of this fibration.

Each fixed point gives rise to a trivial section of the product bundle
$$M^{G}\rightarrow BG\times M^{G}\rightarrow BG.$$
Composing with the injection
$M_{G}\times BG\hookrightarrow EG\times M/G=M_{hG}$
gives a section of the Borel fibration. So we have a natural injection:
$$k: M^{G}\hookrightarrow M^{hG}.$$

For any $G$-CW complex $M$, there is an equivariant rationalization $m: M\rightarrow M_{\mathbb{Q}}$, that is, $M_{\mathbb{Q}}$ is also a $G$-CW complex, $m$ is an equvariant map, and $(M_{\mathbb{Q}})^G\simeq (M^{G})_{\mathbb{Q}}.$ Moreover, we have
\begin{proposition}\cite[Proposition 12]{U-Y-S}
If $M$ is a Postnikov piece, that is $\pi_{>N}(M)=0$ for some N, then

(i) $M^{hG}$ has the homotopy type of a nilpotent CW-complex of finite type.

(ii) $(M^{hG})_{\mathbb{Q}}\simeq (M_{\mathbb{Q}})^{hG}$.
\end{proposition}

Note that if $M_{\mathbb{Q}}$ is a Postnikov piece, then $(M_{\mathbb{Q}})^{hG}$ makes sense and is a rational space.

Now, we determine the homotopy type of the homotopy fixed point sets of certain $S^3$-actions. That is, we give the proof of Theorem \ref{2} and \ref{4}.
\begin{proof}[Proof of Theorem \ref{2}]

(1) \textbf{Case 1}: $n$ is odd.

We only prove the case $n=4k+3$, the case $n=4k+1$ is similar, so we omit it.

As in the proof of \cite[Theorem 19]{U-Y-S}, it is not hard to get the model of the corresponding Borel fiberation
$$\xi:(A,0)\hookrightarrow((\Lambda e)\otimes A,D)\rightarrow (\Lambda e,0),$$
in which $(A,0)=(\Lambda x/x^{k},0)$ and $|x|=4$, $|e|=n$.
This fibration is trivial, so
Sec$(\xi)\simeq{\rm{Map}}(\mathbb{H}P^{k},S^n)$.

By \cite[Theorem 9]{U-Y-S}, the model of ${S^{n}_{\mathbb{Q}}}^{hS^3}$ is $(\Lambda(x_{1},x_{2},\cdots, x_{n+1/4}),0)$. It is exactly the model of $S^3\times S^{7}\times\cdots \times S^n$. It follows that ${S^{n}_{\mathbb{Q}}}^{hS^3}\simeq_{\mathbb{Q}} S^a\times S^{a+4}\times\cdots \times S^n.$

(2) \textbf{Case 2}: $n=4k$.

As $\pi_{\geq 2n}(S^n)\otimes \mathbb{Q}=0$, a model of the Borel fibration is
$$\xi_{2n}:(A,0)\hookrightarrow(\Lambda(e,e')\otimes A,D)\rightarrow (\Lambda(e,e'),d),$$
where $A=\Lambda x/x^{2k+1}$, $x$, $e$, $e'$ are of degree $4$, $n$, $2n-1$ respectively, $De=0$, $De'=e^2+\lambda x^{\frac{n}{4}}e$, $de'=e^2$.

(i) If $\lambda=0$, then $\xi_{2n}$ is trivial and
$${S^{n}_{\mathbb{Q}}}^{hS^3}\simeq {\rm{Map}}(\mathbb{H}P^{2k}, S^n)_{\mathbb{Q}}.$$
A straightforward computation shows that this mapping space has a model of the form
$$(\Lambda y_1,0)\otimes(\Lambda ((x_{s})_{1\leq s\leq k},(y_r)_{2\leq r\leq 2k}),d)$$
with $|x_s|=4s$, $|y_{r}|=4r-1$, $dx_s=0$ $(1\leq s\leq k)$,
$\displaystyle{dy_r=\sum_{s+t=r}x_sx_t}$ $(r>1).$

(ii) If $\lambda\neq 0$, then the fibration $\xi_n$ has two non homotopic sections $\sigma$, $\tau$ which corresponds to the only two possible retractions of its model,
$$\varphi_{\sigma},\varphi_{\tau}:(\Lambda(e,e')\otimes A,D)\rightarrow (A,0),\ \ \ \ \ \ \ \varphi_{\sigma}(e)=0,\ \ \varphi_{\tau}(e)=\lambda x^{k}.$$ By the same way in \cite{U-Y-S}, we have that the model of ${\rm{Sec}}_{\sigma}(\xi_{2n})$ is of the form
$$(\Lambda ((x_{s})_{1\leq s\leq k},(y_r)_{1\leq r\leq 2k}),\widetilde{d})$$
with $|x_s|=4s$, $|y_r|=4r-1$. The linear part of $\widetilde{d}$ is:
$$\widetilde{d}(y_{r})=\lambda x_{r}$$
for $1\leq r\leq k$. Which means the minimal model of ${\rm{Sec}}_{\sigma}(\xi_{2n})$ is
$$(\Lambda(y_r)_{k+1\leq r\leq 2k}, 0).$$

Replace $\lambda$ by $-\lambda$, we have that the model of ${\rm{Sec}}_{\tau}(\xi_{2n})$ is the same.

(3) \textbf{Case 2}: $n=4k+2$.

As $\pi_{\geq 2n}(S^n)\otimes \mathbb{Q}=0$, a model of the Borel fibration is
$$\xi_{2n}:(A,0)\hookrightarrow(\Lambda(e,e')\otimes A,D)\rightarrow (\Lambda(e,e'),d),$$
where $A=\Lambda x/x^{2k+1}$. $x$, $e$, $e'$ are of degree $4$, $n$, $2n-1$ respectively, $De=0$, $De'=e^2$, $de'=e^2$.

So the fibration $\xi_{2n}$ is trivial, we have
$${S^{n}_{\mathbb{Q}}}^{hG}\simeq {\rm{Map}}(\mathbb{H}P^{2k}, S^n)_{\mathbb{Q}}.$$

The model of ${S^{n}_{\mathbb{Q}}}^{hG}$ is:
$$(\Lambda(y_1,y_2),0)\otimes(\Lambda ((x_{s})_{1\leq s\leq k},(y_r)_{3\leq r\leq 2k+1}),d)$$
with $|x_s|=4s+2$, $|y_{r}|=4r-1$, $dx_s=0$ $(1\leq s\leq k)$,
$\displaystyle{dy_r=\sum_{s+t=r-1}x_sx_t}$ $(3\leq r\leq 2k+1).$

This implies the result.
\end{proof}
\begin{proof}[Proof of Theorem \ref{4}]
First, we assume $n=2k+1$. As $\pi_{\geq4k+4}(\mathbb{C}P^n_{\mathbb{Q}})=0$, we just need to use the model of $\xi_{2n+2}$:
$$(A,0)\rightarrow(\Lambda(e,e')\otimes A, D)\rightarrow(\Lambda (e,e'),d),$$
in which $A=(\Lambda x)/x^{k+2}$, $|x|=4$, $|e|=2$, $|e'|=4k+3$, and
$$De=0,\ \ \ \ \   De'=e^{n+1}+\sum_{j=1}^{k}\lambda_{j}e^{j}x^{n+1-2j},\ \ \ \ \lambda\in\mathbb{Q}, j=1,\cdots,n.$$
The retraction of this model of fibration is just $\varphi(e)=0$. So we have ${\rm{Sec}}(\xi_{4k+4})$ is connected, and the model of it is
$$(\Lambda(e,(e'_{r})_{1\leq r\leq k+1}, \widetilde{d})$$
with $|e|=2$, $|e'_r|=4r-1$, $\widetilde{d}(e'_r)=\lambda_{k+1-r} e^{2r}$ for $1\leq r\leq k$ and $\widetilde{d}(e'_{k+1})=e^{2k+2}$.

If $\lambda_1\neq 0$ this is a model of
$$S^2\times S^7\times\cdots \times S^{4k+3}.$$
If $\lambda_1=\cdots=\lambda_{i-1}=0$ and $\lambda_i\neq 0$, this is a model of
$$S^3\times\cdots \times S^{4k-4i-1} \times CP^{2k+1-2i}\times S^{4k-4i+3}\times \cdots \times S^{4k+3}.$$
Finally, if all $\lambda_i=0$, then it is a model of
$$S^3\times S^{7} \times \cdots\times S^{4k-1} \times CP^{2k+1}.$$

For $n$ even, the proof is similar, so we omit it.
\end{proof}
\section{The Inclusion $k:M^{S^1}\hookrightarrow M^{hS^1}$}
We begin with some interesting observations on $S^1$-actions.

In \cite[Example 12]{U-Y-A}, there is an $S^1$-action on $M=K(\mathbb{Z},n)\times K(\mathbb{Z}, n+1)$, such that the model of it's Borel fibration is
$$\eta_{n}:(\Lambda x,0)\hookrightarrow (\Lambda x\otimes\Lambda(z,y),D)\rightarrow(\Lambda (z,y),d),$$
where $|x|=2$, $|z|=n$, $|y|=n+1$, $D(z)=0$, and $D(y)=xz$. For $n=2k$, there is only one retraction $\sigma$: $\sigma(z)=\sigma(y)=0$, so ${\rm{Sec}}(\eta_{2k})$ is path connected.

By the same method used in \cite{U-Y-S}, a model of ${\rm{Sec}}(\eta_{2k})$ is
$$(\Lambda ((z_i)_{1\leq i\leq k}, (y_j)_{1\leq j\leq k+1}),d),$$
where $|z_i|=2i$, $|y_j|=2j-1$ and $d(y_i)=z_{i}$. The minimal model of ${\rm{Sec}}(\eta_{2k})$ is $(\Lambda y_{k+1},0)$, so ${\rm{Sec}}(\eta_{2k})\simeq_{\mathbb{Q}}S^{2k+1}$ is an elliptic space. However, $M$ is not an elliptic space.

Next we complement \cite[Corollary 2]{U-Y-S} with
the following:
\begin{proposition}\label{1}
For an $S^1$-space $M$ which is a nilpotent finite complex, the following conditions are equivalent:\\
1). $M$ is elliptic, \\
2). Each component of ${M_{\mathbb{Q}}}^{hS^1}$ is elliptic,\\
3). One of the components of ${M_{\mathbb{Q}}}^{hS^1}$ is elliptic.
\end{proposition}
\begin{proof}[Proof of Proposition \ref{1}]
$\\$
$1)\Rightarrow 2)$: \cite[Theorem 15]{U-Y-S}.\\
$2)\Rightarrow 3)$: Trivial.\\
$3)\Rightarrow 1)$: By \cite[Theorem 13]{U-Y-A}, $2{\rm{dim}}\ \pi_{\ast}({\rm{Sec}}_{\sigma}(\xi)\otimes \mathbb{Q})\geq {\rm{dim}}\ \pi_{\ast}(M)\otimes \mathbb{Q}.$
By ${\rm{Sec}}_{\sigma}(\xi)$ is elliptic, ${\rm{dim}}\ \pi_*({\rm{Sec}}_{\sigma}(\xi))\otimes \mathbb{Q}$ is finite, so ${\rm{dim}}\ \pi_*(M)\otimes \mathbb{Q}$ is finite. Then $M$ is elliptic.
\end{proof}
\begin{remark}
\emph{The theorem holds also for $G=S^3$. The proof is similar.}
\end{remark}

The rest of the section is devoted to the proof of Theorem \ref{5}.

Let $M$ be an $S^1$-space and $M^G\neq\emptyset$. Then the inclusion $M^{S^1}\hookrightarrow M$ induces a map of Borel fibrations:

\begin{equation}\label{D1}
\xymatrix{
M^{S^1}\ar[rr]\ar[d]&&M\ar[d]\\
\mathbb{C}P^\infty\times M^{S^1}\ar[rr]^{\gamma}\ar[rd]_{\eta}&&M_{hS^1}\ar[ld]^{\xi}\\
&\mathbb{C}P^\infty.&
}
\end{equation}
If there exists $N$ such that $\pi_{\geq N}(M_{\mathbb{Q}})=0$ and $\pi_{\geq N}(M^{S^1}_{\mathbb{Q}})=0$. Then $k$ is identified with the corresponding
$$M^{S^1}\hookrightarrow{\rm{Map}}((\mathbb{C}P^{\infty})^{(N)}, M^{S^1})\rightarrow {\rm{Sec}}(\xi_{N})\cong M^{hS^1},$$
obtained by truncating in the diagram (\ref{D1}):

\centerline{
\xymatrix{
M^{S^1}\ar[rr]\ar[d]&&M\ar[d]\\
F_N\ar[rr]^{\gamma_N}\ar[rd]_{\eta_N}&&E_N\ar[ld]^{\xi_N}\\
&(\mathbb{C}P^\infty)^{(N)}.&
}
}
Now let
\begin{equation}
\xymatrix{
&(A\otimes \Lambda V,D)\ar[r]\ar[dd]_{\psi}&(\Lambda V,d)\ar[dd]_{\varphi}\\
(A,0)\ar[ur]\ar[dr]&&\\
&(A,0)\otimes (\Lambda Z,d)\ar[r]& (\Lambda Z,d)
}
\end{equation}
be a model of the above diagram where $(A,0)=(\Lambda x/ (\Lambda x)^{>N},0)$, $(\Lambda V,d)$ and $(\Lambda Z,d)$ are minimal Sullivan models of $M$ and $M^{S^1}$, respectively.

Then we have the following theorem:
\begin{theorem}\cite[Theorem 21]{U-Y-S}\label{t}
The composition

\centerline{
\xymatrix{
(\Lambda(V\otimes A^{\#}),\widetilde{d})\ar[r]^{\phi}&(\Lambda(Z\otimes A^{\#}),\widetilde{d})\ar[r]^{\gamma}&(\Lambda Z,d)
}
}
is a model of $k:M_{\mathbb{Q}}^{S^1}\hookrightarrow M_{\mathbb{Q}}^{hS^1}$. And the morphism are defined by:
$$\phi(v\otimes \alpha)=\rho^{-1}[\psi(v)\otimes \alpha],\ \ v\otimes\alpha\in V\otimes A^{\#},$$
$$\gamma(z\otimes \alpha)=\left\{
\begin{array}{cc}
z&  \alpha=1,\\
0&   \alpha\neq 1,
\end{array}
\right.
\ \ \ \ z\otimes\alpha\in Z\otimes A^{\#}.$$
\end{theorem}
Then we give some information about $\psi$. First,  let $(\Lambda x\otimes \Lambda V,D)$ be a model of the fibration $\xi$, we can decompose the differential $D$ in $A\otimes \Lambda V$ into
$$D=\sum_{i\leq 1}D_i,\ \ \ \ D_{i}(V)\subset \Lambda x\otimes \Lambda^i V.$$

And we have:
\begin{proposition}\cite[Lemma 14]{U-Y-A}\label{l1}
The vector space $V$ can be decomposed into a direct sum $W\oplus K\oplus S$ where \\
(1) $W\oplus K=ker D_1$,\\
(2) $K$ and $S$ have the same dimension admitting bases $\{v_{i}\}_{i\in I}$, $\{s_{i}\}_{i\in I}$, and for any $i\in I$, there exists $n_{i}\geq 1$ such that $D_{1}(s_{i})=x^{n_{i}}v_{i}$.
\end{proposition}
Let $\mathbb{K}=Q(x)$, the field of fractions of $\Lambda x$, we obtain a morphism of (ungraded) differential vector spaces
$$\overline{\psi}:(\mathbb{K}\otimes V, D_1)\rightarrow (\mathbb{K}\otimes Z, 0)=(Z_{\mathbb{K}},0).$$
If we assume $\mathbb{K}$ concentrated in degree 0 and consider in $V$  and $Z$ the usual $\mathbb{Z}_2$-grading given by the parity of the generators, then the Borel localization theorem claim that:
\begin{theorem}\cite[Theorem 22]{U-Y-S}
The morphism
$$\overline{\psi}:(\mathbb{K}\otimes V, D_1)\rightarrow (Z_{\mathbb{K}},0)$$
is a quasi-isomorphism.
\end{theorem}
By Proposition \ref{l1}, we have
\begin{lemma}\label{l2}
(1). {\rm{dim}} $W$={\rm{dim}} $Z$.

(2). There are $\{w_{j}\}_{j\in J}$, $\{z_{j}\}_{j\in J}$ which are homogenous basis of $W$ and $Z$ respectively, and non negative integers $\{m_j\}_{j\in J}$ such that
$$\psi(w_j)=x^{m_{j}}z_j+\Gamma_j,\ \ \ \ \Gamma_j\in R \otimes \Lambda^{\geq2}Z,\ \ \ \ j\in J,$$
and
$$\psi(s_{i})\in R\otimes\Lambda^{\geq2}Z,\ \ \ \ \ \psi(v_{i})\in R\otimes\Lambda^{\geq2}Z,\ \ \ \ s_{i}\in S,\  v_{i}\in K,\ i\in I. $$
\end{lemma}
\begin{theorem}
For an $S^1$-complex $M$ which is simply connected with
$${\rm{dim}}\ \pi_*(M)\otimes \mathbb{Q}<\infty.$$
Then the inclusion
$$k:M^{S^1}\hookrightarrow M^{hS^1}$$
is a rational homotopy equivalence if and only if $M$ is rational homotopy equivalent to a product of $CP^{\infty}$.
\end{theorem}
\begin{proof}
By Theorem \ref{t}, the model of $k$ is:
$$\alpha:(\Lambda(V\otimes A^{\#}),\widetilde{d})\rightarrow(\Lambda(Z\otimes A^{\#}),\widetilde{d})\rightarrow(\Lambda Z,d).$$
By \cite[Theorem 24]{U-Y-S}, $\pi_{*}(k)\otimes \mathbb{Q}$ is injective, so we only consider the surjective part.

By \cite[Theorem11]{U-Y-S}, $(\Lambda(V\otimes A^{\#}),\widetilde{d})$ is a model of $M_{\mathbb{Q}}^{hS^1}$. Then we have
$$H^k(V\otimes A^{\#}, \widetilde{d}_{1})\cong Hom(\pi_{k}(M_{\mathbb{Q}}^{hS^1}),\mathbb{Q}),$$
where $k\geq 1$.

By Proposition \ref{l1}, $V=W\oplus K\oplus S$. An easy computation shows that $(W\otimes A^{\#})\oplus S\subset H^*(V\otimes A^{\#}, \widetilde{d}_{1})$.
It is obvious that
$$\alpha(w_{j})=0\Leftrightarrow m_{j}\neq 0,$$
$$\alpha(w_j\otimes (x^i)^{\#})=0\Leftrightarrow m_{j}\neq i,$$
$$\alpha(s_j)=0.$$

If there exists $j$ such that $|w_{j}|\geq 2$ or $S\neq\emptyset$, then $H(\alpha,\widetilde{d}_{1})$ is not injective, so $k$ is not a rational homotopy equivalence.

If $|w_{j}|=2$, for each $j\in J$, and $S=\emptyset$, we have $(\Lambda W,d)$ is a model of a product of $\mathbb{C}P^{\infty}$. It is easy to show that $k$ is a rational homotopy equivalence.
\end{proof}


\begin{thebibliography}{10}
\bibitem{U-Y-S}U. Buijs, Y. F\'{e}lix, S. Huerta, A. Murillo, {\em The homotopy fixed point set of a Lie group actions on elliptic spaces}, Proc. London. Math. Soc. \textbf{110} (5), 1135-1156, 2014.
\bibitem{U-Y-A}U. Buijs, Y. F\'{e}lix, A. Murillo, {\em Rational homotopy of the (homotopy) fixed point sets of circle actions}, Adv. Math. \textbf{222} (1), 151-171, 2009.
\bibitem{U-A}U. Buijs, A. Murillo, {\em Basic constructions in rational homotopy theory of function spaces}, Ann. Inst. Fourier \textbf{56} (3), 815-838, 2006.
\bibitem{Y-S}Y. F\'{e}lix, S. Halperin, J.-C. Thomas, {\em Rational homotopy theory}, Graduate Texts in Mathematics \textbf{205}, Springer-Verlag, New York, 2001.
\bibitem{J}J. Goyo, {\em The Sullivan model of the homoyopy-fixed-point set}, thesis, University of Toronto, 1989.
\bibitem{A}A. Heafliger, {\em Rational honotopy of the space of sections of a nilpotent bundle}, Trans. Amer. Math. Soc. \textbf{273} (2), 609-620, 1982.
\bibitem{H}H. Miller, {\em The Sullivan conjecture on maps from classifying spaces}, Ann. of Math. \textbf{120} (1), 39-87, 1984.
\end{thebibliography}
\end{document}